\documentclass[12pt,a4paper]{article}
\usepackage{amsmath,amssymb,amstext,amsfonts,latexsym,amsthm}
\usepackage{dsfont,mathrsfs,txfonts}
\usepackage{graphicx,hyperref,color}
\usepackage[affil-sl]{authblk}

\newtheorem{theorem}{Theorem}[section]
\newtheorem{lemma}{Lemma}[section]
\newtheorem{proposition}{Proposition}[section]
\newtheorem{definition}{Definition}[section]
\newtheorem{corollary}{Corollary}[section]

\newtheorem{example}{Example}[section]

\newcounter{Th-Alfa}
\newcommand{\XX}{\mathds{X}}
\newcommand{\CC}{\mathds{C}}
\newcommand{\MM}{\mathds{M}_n}
\newcommand{\NN}{\mathds{N}}
\newcommand{\PP}{\mathds{P}}

\newcommand{\RR}{\mathds{R}}
\newcommand{\ZZ}{\mathds{Z}}
\newcommand{\ZZp}{\mathds{Z}_+}
\newcommand{\dsty}{\displaystyle}
\newcommand{\sgn}[1]{{\mathop{\rm sgn}\!}\left(#1 \right)}
\newcommand{\supp}{\mathop{\rm supp}}

\newcommand{\dgr}[1]{\mbox{{\rm deg\/}}(#1)}

\newcommand{\inter}[1]{#1\strut^{\mathrm{o}}\!}
\newcommand{\intersub}[1]{#1\strut^{\!\!\!\mathrm{o}}}
\newcommand{\ch}[1]{\mathbf{Co}\!\left(#1 \right)} 
\newcommand{\intch}[1]{\ch{#1}\strut^{\!\!\mathrm{o}}} 

\DeclareMathOperator*\lowlim{\underline{lim}}
\DeclareMathOperator*\uplim{\overline{lim}}
\newcommand{\cp}[1]{\mathrm{cap}\!\left(#1 \right)}
\DeclareRobustCommand{\wlim}{\mathop{\operatorname{w{-}lim}}}
\def\bbuildrel#1_#2^#3{\mathrel{
 \mathop{\kern 0pt#1}\limits_{#2}^{#3}}}
\def\bbbuildrel#1_#2{\mathrel{
 \mathop{\kern 0pt#1}\limits_{#2}}}
\newcommand{\MinDeg}[1][M]{\mathfrak{u}_{#1}}
\newcommand{\Nzero}[2]{\mathbf{N}_z\!\left(#1;#2\right)}
\newcommand{\Ncero}[2]{\mathbf{N}_{o}\!\left(#1;#2\right)}

\newcommand{\funD}[4]{\begin{cases} #1 , & \hbox{if } \; #2;\\
										#3, & \hbox{if } \; #4 \end{cases}}
\def\bbuildrel#1_#2^#3{\mathrel{
 \mathop{\kern 0pt#1}\limits_{#2}^{#3}}}
\def\bbbuildrel#1_#2{\mathrel{
 \mathop{\kern 0pt#1}\limits_{#2}}}

\title{\textbf{Polynomials of least deviation from zero in  Sobolev $p$-norm}}
\author[1]{Abel D\'{\i}az-Gonz\'{a}lez\thanks{abdiazgo@math.uc3m.es}}
\author[1]{H\'{e}ctor Pijeira-Cabrera\thanks{hpijeira@math.uc3m.es}\thanks{Research partially supported by  Ministry of Science, Innovation and Universities of Spain, under grant  PGC2018-096504-B-C33}}
\author[1]{Javier Quintero-Roba\thanks{jaquinte@math.uc3m.es}} 
\affil[1]{Universidad Carlos III de Madrid, Spain}

\date{}
\begin{document}

\maketitle

\begin{abstract}
The first part  of this paper  complements previous  results on  characterization of polynomials of least deviation from zero in  Sobolev $p$-norm ($1<p<\infty$) for the case $p=1$. Some relevant  examples are indicated.

The second part deals with  the  location of zeros of polynomials of least deviation  in  discrete Sobolev $p$-norm. The  asymptotic distribution of zeros is established on  general conditions. Under some  order restriction  in the discrete part, we prove   that,  the $n$-th polynomial of least deviation   has at least $n-\mathbf{d}^*$ zeros on the convex hull of the support of the measure, where $\mathbf{d}^*$ denotes the number of terms in the discrete part.

\medskip
\begin{description}
  \item[Keywords:] polynomials of least deviation from zero, extremal polynomials, Sobolev norm, zero location
  \item[2020 MSC:] 30C10, 30C15, 33C47, 41A10, 41A50, 41A52, 46E35
\end{description}
\medskip

\end{abstract}

\section{Introduction}

Let $\PP$ be the linear space of polynomials,  $\|\cdot\|$ be a norm defined on $\PP$ and $\PP^1_n$ be  the subset of all polynomials of degree $n \in \ZZp$ whose  leading coefficient is equal to one (monic). A classic problem in analysis is the existence, uniqueness and characterization of the monic polynomial of degree $n \in \ZZp$ with minimum deviation from  zero with respect to the norm  $\|\cdot\|$, i.e. the polynomial $P_n(z)=z^n+ \dots$ such that
\begin{equation}\label{DefExtPoly}
  \|P_n\|= \inf_{Q_n \in \PP^1_n }  \|Q_n\|.
\end{equation}
A polynomial $P_n \in \PP^1_n$  that satisfies \eqref{DefExtPoly}  is called   polynomial of  least deviation from zero with respect to $\|\cdot\|$,   for brevity, a  $n$-th minimal (or extremal) polynomial with respect to  $\|\cdot\|$. This problem has its origin in the study carried out by P. L. Chebyshev on the decrease of the  friction in the joints of the Watt parallelogram that converts the movement of the piston of the steam engine into wheel rotation. As a consequence, what we know today as Chebyshev polynomials were discovered (c.f.  \cite[Ch. 1]{Bog13}). It is well known that Chebyshev monic  polynomials of the first kind  are minimal with respect to the uniform norm at $ [- 1,1] $ and that those of the second kind are minimal with respect to the usual norm at $L^1[-1,1]$ (c.f. \cite[\S6.6]{Cheney66} or \cite[\S 3.3]{Dav75} ). Let us mention that these works constituted a starting point   of the general theory of orthogonal polynomials. Today,    minimal polynomials are of great interest in various areas such as approximation theory, potential theory, optimization of numerical algorithms, and signal processing.

Note that, any polynomial $Q \in \PP^1_n$ could be  written as $Q(z)=z^n-q(z)$ with  $q \in \PP_{n-1}$. Let $q_0$ be a fixed element of $\PP_{n-1}$ and define the associated subset

$$\mathds{A}_{n,0}=\{q \in \PP_{n-1}: \|x^n-q\| \leq \|x^n-q_0\|\}.$$

As $\mathds{A}_{n,0}$ is  a compact subset of $\PP_{n-1}$, there exists $q_1 \in \PP_n$ such that $\|x^n-q_1\| \leq \|x^n-q\|$ for all  $q \in \PP_{n-1}$, in virtue of the arbitrariness of $q_0$. Hence, the existence of a minimal polynomial  is  guaranteed. However, the uniqueness of the minimal polynomial with respect to \eqref{SobNorm01} is not always ensured, as we will show in some of our case studies.

Nevertheless, it is straightforward to prove that $\MM$ (the set of all monic  minimal polynomials with respect to  $\|\cdot\|$ of degree $n$)  is a convex set. Indeed, if $Q_n, R_n \in \MM$ and $\lambda \in [0,1]$,  then $P_n(x)=\lambda Q_n+ (1-\lambda) R_n(x)$ is also an element of $\MM$ since

$$
\|P_n\|=\|\lambda Q_n+ (1-\lambda) R_n(x)\|  \leq  \lambda\|Q_n\|+(1-\lambda)\|R_n\|=\|Q_n\|.
$$

In this paper, we are interested in the case in which the norm $\|\cdot\|$    is as we define below. Let $1\leq p<\infty$ and consider the vector of measures $\vec{\mu}= (\mu_0,\mu_1,\dots,\mu_m)$, for $m\in\ZZp$, where $\mu_k$ is a positive finite Borel measure with $\supp\mu_k\subset \RR$ and $\PP\subset L^1\left(\mu_k\right)$ for $k=0,1,\dots,m$. Denote by $\Delta_k$ the convex hull of $\supp\mu_k$, that is the smallest  interval containing $\supp\mu_k$.  Let  $f^{(k)}$ denote  the $k$-th derivative of a function $f$. If $\Delta_0$ contains infinite elements, the expression
\begin{equation}\label{SobNorm01}
\|f\|_{p,\vec{\mu}}=  \left(\sum_{k=0}^{m} \|f^{(k)} \|_{k,p}^p\right)^{1/p} =\left(\sum_{k=0}^m \int_{\Delta_k} \left|f^{(k)}\right|^p d\mu_k\right)^{1/p},
\end{equation}
defines a norm over $\PP$ known as the Sobolev $p$-norm and the vector of measures $\vec{\mu}$ is called standard. If each  measure $\mu_k$,  $0\leq k \leq m$ satisfies  $\mu_k(\{x\})=0$ for all $x \in \RR$, we say that the vector of measures $\vec{\mu}$ is continuous.

First, observe that  for $m=0$ this norm reduces to the usual $L^p\left(\mu_0\right)$ norm. We will call \emph{$n$-th Sobolev minimal polynomial with respect to $\|\cdot\|_{p,\vec{\mu}}$}, to any polynomial $P_n \in \PP^1_n$ that is a solution of the minimal problem \eqref{DefExtPoly}.

For the norm \eqref{SobNorm01} with $\vec{\mu}$ standard, we consider  two different cases:
  \begin{description}
     \item[\textbf{Continuous Sobolev norms,}]  if  $\vec{\mu}$ is continuous.
 \item[\textbf{Discrete Sobolev norms,}] if for every $k=1, \ldots, m$ the measure $\mu_k$ is supported on a finite number of points.
\end{description}

It is said that a   Sobolev $p$-norm   is  \emph{sequentially dominated}  if $\supp\mu_k\subset \supp\mu_{k-1}$
and $ d\mu_k = f_{k-1}d\mu_{k-1}$  where $f_{k-1} \in L_{\infty}(\mu_{k-1})$ and $k = 1,\ldots,m$. Furthermore, the norm  \eqref {SobNorm01} on $\PP$  is said to be \emph{essentially sequentially dominated},  if there exists    a sequentially dominated norm that is equivalent to \eqref {SobNorm01}. As usual, two norms $\|\cdot\|_1$ and $\|\cdot\|_2$ on a given normed space $\XX$ are said to be equivalent if there exist positive constants $c_1,c_2$ such that
$ c_1\|x\| \leq \|x\| \leq c_2\|x\|$ for all $x\in \XX$.

The notions of \emph{sequentially dominated norm} and \emph{essentially sequentially dominated norm} were  introduced  in \cite{HPLago99} and \cite{Jomaro01} respectively. Both   notions are closely related to the uniform boundedness of the distance between  the zeros of sequences of minimal polynomials and the support of the measures involved in  \eqref{SobNorm01}.  For more details on this aspect in the continuous case, we refer the reader to \cite{DurSaff01,LagIpHp01} for  $p=2$, \cite{GlIpHp05,Jomaro12} for $1<p<\infty$ and \cite{DiOrPi08,DiOrPi09,LopMarPij06}   for $p=2$ and measures with unbounded support.

Let  $N \in \ZZp$, $\Omega=\{ c_1,\dots,c_N\}\subset \CC $, $\{m_0,\dots,m_N\} \subset \ZZp$ and $m=\max\{m_0,\dots,m_N\}$.  In the discrete case, we will restrict our attention to  Sobolev $p$-norm  under the  following assumptions:
\begin{itemize}
  \item $\dsty \mu_0=\mu+\sum_{j=1}^N  A_{j,0}\delta_{c_j}$, where   $ A_{j,0}\geq 0$, $\mu$ is  a finite positive Borel measure,  $\supp{\mu} \subset \RR$  with infinitely many points,  $\PP\subset L^{1}(\mu)$ and $\delta_{x}$ denotes the Dirac measure with mass one at the point $x$.
  \item For $k=1,\dots,m$;  $\dsty \mu_k=\sum_{j=1}^N A_{j,k}\delta_{c_j}$ where  $A_{j,k}\geq  0$, $A_{j,m_j}>0$,   and  $A_{j,k}=0$ if  $m_j<k\leq m$.
\end{itemize}

We say that a discrete Sobolev $p$-norm   is \emph{non-lacunary} if  $A_{j,k}>0$ for all $0\leq k \leq m_j$ and  $0 \leq j \leq N$.   In any other case, we say that the  discrete Sobolev $p$-norm    is \emph{lacunary}. Obviously, a discrete Sobolev $p$-norm is non-lacunary if and only if is sequentially dominated. A discrete Sobolev $p$-norm  is \emph{essentially non-lacunary} if it is equivalent to a non-lacunary norm.

It is known that the minimal polynomial in $L^p(\mu_0)$ spaces ($m=0$) satisfies the following characterization (see \cite[Sec.2.2 and Ex 7-h]{Borwein95}). A monic polynomial $P_n$ is the $n$-th minimal polynomial in $L_p(\mu_0)$ if and only if

\begin{align*}
 \langle P_n , q\rangle_{p,\mu_0} = &  \int_{\Delta_{0}} q \operatorname{sgn}\left(P_{n}\right)\left|P_{n}\right|^{p-1} d \mu_{0}=0\; \text{for  all }q\in\PP_{n-1},  \\
 \text { where } & \operatorname{sgn}(y)=\left\{\begin{array}{ll}
y /|y|, & \text { if } y \neq 0 ; \\
0, & \text { if } y=0.
\end{array}\right.
\end{align*}

In \cite[Th.4]{AbelPijeiraLago17}, the authors provide the following  extension of this characterization to the Sobolev case when $1<p<\infty$.

\begin{theorem}\label{TheoremPProduct} Consider the Sobolev $p$-norm \eqref{SobNorm01} for $1<p<\infty$. Then the monic polynomial $P_n$ is the $n$-th Sobolev minimal polynomial  if and only if
\begin{equation}\label{CNS-1<p<Inf}
\langle P_n , q\rangle_{p,\vec{\mu}}=\sum_{k=0}^{m} \int_{\Delta_k} q^{(k)} \operatorname{sgn}\left(P_{n}^{(k)}\right)\left|P_{n}^{(k)}\right|^{p-1} d \mu_{k}=0,
\end{equation}
for every polynomial $q\in\PP_{n-1}$.
\end{theorem}

The results in this work complement previous ones in \cite[\S 2]{AbelPijeiraLago17}. There, for $ 1 <p <\infty $,  Theorem \ref{TheoremPProduct}, Proposition \ref{PropUnicidad} and Corollary \ref{CoroZeros} were proved.

The aim of  Section \ref{Sect2}  is to extend   Theorem \ref{TheoremPProduct} to  the case $p=1$. In Theorem \ref{ThChar-p}, we give a general sufficient condition for existence of a minimal polynomial with respect to \eqref{SobNorm01} ($ 1\leq p <\infty $). For $p=1$ this condition is not necessary, as we show in  Examples \ref{ExamNoUniqDisc} and \ref{UniciNoCond}. Furthermore, Example \ref{ExamNoUniqCont} shows that it does not guarantee uniqueness either.  Theorem \ref{characterizationP1} establishes a necessary and sufficient condition  under which \eqref{CNS-1<p<Inf} characterizes minimality with respect to \eqref{SobNorm01} when   $p=1$.

The last two sections deal with  discrete Sobolev norms.  In Section \ref{Sec3}, for essentially non-lacunary Sobolev norms, we give a sufficient condition for the uniform boundedness of the set of zeros of a sequence on minimal polynomials $\{P_n\}$ (see Theorem \ref{ThEssNonL}). Moreover, the asymptotic distribution of zeros is established in Theorem \ref{AsymZerosDis}. Finally, in Section \ref{Sec4},  we introduce the notion of sequentially-ordered Sobolev $p$-norm.  Under this assumption, we prove Theorem \ref{ThZeroLocDiscrete},   which generalizes several known results on the number of zeros of  the $n$-th polynomial of least deviation inside  the convex hull of the support of the measure $\mu$.

\section{Polynomials of least deviation from zero when  $p=1$} \label{Sect2}

Let us first recall a  basic property of the Sobolev norm \eqref{SobNorm01}. Let $R$ be a monic  polynomial with complex coefficients, and let  us write $R=R_1+iR_2$, where  $R_1$ and $R_2$ are polynomials with real coefficients.  Note that $R_1$ is also a monic polynomial with the same degree of $R$ and satisfying
\begin{align*}
\|R\|_{p,\vec{\mu}} ^p&=\sum_{k=0}^{m}\int_{\Delta_k} |R_1^{(k)}+iR_2^{(k)}|^p d\mu_k=\sum_{k=0}^{m}\int_{\Delta_k} \left(\left(R_1^{(k)}\right)^2+\left(R_2^{(k)}\right)^2\right)^{p/2}d\mu_k\\&> \sum_{k=0}^{m}\int_{\Delta_k}\left| R_1^{(k)}\right|^pd\mu_k=\|R_1\|_{p,\vec{\mu}} ^p.
\end{align*}
Therefore,  any  $n$-th Sobolev minimal polynomial with respect to $\|\cdot\|_{p,\vec{\mu}}$, has real coefficients.

 \begin{proposition}[{\cite[Prop. 1]{AbelPijeiraLago17}}] \label{PropUnicidad} Let $\|\cdot\|_{p,\vec{\mu}}$  be the Sobolev type norm defined by \eqref{SobNorm01}, with $1<p< \infty$. Then,  there exists a unique $P_{n} \in \PP^1_n$  such that $\dsty \|P_{n}\|_{p,\vec{\mu}}= \inf_{Q_n \in\PP^1_{n}} \|Q_n\|_{p,\vec{\mu}}.$
\end{proposition}

\begin{theorem}[Sufficient condition]\label{ThChar-p}
Consider the Sobolev $p$-norm \eqref{SobNorm01} for $1\leq p<\infty$, when  $\vec{\mu}=(\mu_0,\dots,\mu_m)$ is a standard vector measure. If $P_n \in \PP^1_n$ is  such that for all $q\in \PP_{n-1}$
\begin{equation}\label{sufcond}
\langle P_n , q\rangle_{p,\vec{\mu}}=\sum_{k=0}^{m} \int_{\Delta_k} q^{(k)}(x) \, \sgn{P_{n}^{(k)}(x)}\left|P_{n}^{(k)}(x)\right|^{p-1} d \mu_{k}(x)=0,
\end{equation}
then $P_n$ is a  minimal polynomial with respect to $\|\cdot\|_{p,\vec{\mu}}$.
\end{theorem}

\begin{proof} If $1<p<\infty$ the proof is carried out as the proof of the sufficiency   in \cite[Th. 4]{AbelPijeiraLago17}, step by step.

Hence, in what follows we consider $p = 1$. Write $P_n(z)= z^n-q_0(z)$ where  $q_0\in \PP_{n-1}$,   let $q\in\PP_{n-1}$ arbitrary and assume that \eqref{sufcond} holds, then
\begin{align*}
\left\|P_n\right\|_{1,\vec{\mu}} = &\sum_{k=0}^m \int_{\Delta_k} \left(\left(x^n\right)^{(k)}-q_0^{(k)}(x)\right)\sgn{P_n^{(k)}(x)}d\mu_k(x)=\langle P_n , x^n-q_0\rangle_{1,\vec{\mu}}\\
=&\langle P_n , x^n-q+q-q_0\rangle_{1,\vec{\mu}}=\langle P_n , x^n-q\rangle_{1,\vec{\mu}}+\langle P_n , q-q_0\rangle_{1,\vec{\mu}} = \langle P_n , x^n-q\rangle_{1,\vec{\mu}}
\end{align*}
and taking absolute value we have
\begin{align*}
\left\|P_n\right\|_{1,\vec{\mu}} \leq \sum_{k=0}^m \int_{\Delta_k}\left| \left(x^n-q\right)^{(k)}\right| d\mu_k = \|x^n-q\|_{1,\vec{\mu}}, \quad \forall q\in \PP_{n-1},
\end{align*}
which is equivalent to the assertion of the theorem for $p=1$. \qed
\end{proof}

In  \cite[Th. 4]{AbelPijeiraLago17},  it was proved that if $1 <p <\infty$ the condition   \eqref{sufcond} is also necessary, i.e. Theorem \ref{ThChar-p} is a characterization of the extremality  in this case.

With the same arguments as in \cite[Cor. 1 and Cor. 2]{AbelPijeiraLago17},  we have the following corollary.

\begin{corollary} \label{CoroZeros} Under the assumptions of Theorem  \ref{ThChar-p},  if $P_n \in \PP^1_n$  satisfies the condition \eqref{sufcond}, then
\begin{enumerate}
  \item For all  $n\geq 1$,   $P_{n}$  has at least one zero of odd multiplicity  on $\intch{\supp{\mu_0}}$.
  \item  For all  $n\geq 2$,  $P_{n}^{\prime}$   has at least one zero of odd multiplicity  on $\intch{\supp{\mu_0}\cup \supp{\mu_1}}$.
\end{enumerate}
where $\ch{A}$ and $\inter{A}$ denote the  convex hull and the interior of a set $A$, respectively.
\end{corollary}

Observe that   if $p=1$, the condition \eqref{sufcond} only depends on the sign of $P_n$ and its derivatives on the support of the corresponding measure and not on the values of the polynomial itself.  Consequently, unlike what happens in the case $1<p<\infty$, if $p=1$ we lose the uniqueness of the minimal polynomial, as can be seen in the following examples. Furthermore, in Example \ref{ExamNoUniqDisc}, we obtain a minimal  polynomial that does not satisfy the condition \eqref{sufcond}.

\begin{example}[Continuous case]\label{ExamNoUniqCont}\

Consider the Sobolev norm associated to the vector of measures $\vec{\mu}=(\nu|_{[-2,0]},\nu|_{[0,1]})$, where $\nu|_{[a,b]}$ denotes the Lebesgue measure over the real interval $[a.b]$,
\begin{equation}\label{Example-1}
 \|f\|_{1,\vec{\mu}}=\int_{-2}^0|f|dx+\int_{0}^1|f'|dx.
\end{equation}
Let  $ P_{a,2}(x)=(x+1)(x-a)$, with $a \in [ 0,1]$, a family of monic polynomials of degree 2. Note that
\begin{align*}
\langle P_{a,2} , 1\rangle_{1,\vec{\mu}}=&\int_{-2}^0\sgn{(x+1)(x-a)}dx=\int_{-2}^{-1}dx-\int_{-1}^{0}dx=0.\\
\langle P_{a,2} , x\rangle_{1,\vec{\mu}}=&\int_{-2}^0x \, \sgn{(x+1)(x-a)} dx+\int_{0}^1\sgn{2x+1-a}dx \\ =& \int_{-2}^{-1}x\,dx-\int_{-1}^0x\,dx+\int_{0}^1dx=0.
\end{align*}
Then, from Theorem \ref{ThChar-p}, the polynomials $ P_{a,2}$ with $0 \leq a \leq 1$ are all minimal with respect to \eqref{Example-1}.

Furthermore, note that the   minimal polynomials  $P_{a,2}(x)=(x+1)(x-a)$ for all $0\leq a \leq 1$, are the convex combinations of the   minimal polynomials $x^2-1$ and $x^2+x$.
\end{example}

\begin{example}[Discrete case] \label{ExamNoUniqDisc}\

Consider the Sobolev norm associated to $\vec{\mu}=(\nu|_{[-2,0]},\delta_{0})$, where $\delta_{0}$ is  the Dirac measure with mass one at  $x=0$,
\begin{equation}\label{Example-2}
\|f\|_{1,\vec{\mu}}=\int_{-2}^0|f|dx+|f'(0)|.
\end{equation}
Let  $ P_{b,2}(x)=(x+1)(x-b)$, with $b \in [ 0,1)$, a family of monic polynomials of degree 2. Note that
\begin{align*}
\langle  P_{b,2} , 1\rangle_{1,\vec{\mu}}=&\int_{-2}^0\sgn{(x+1)(x-b)}dx=\int_{-2}^{-1}dx-\int_{-1}^{0}dx=0.\\
\langle  P_{b,2} , x\rangle_{1,\vec{\mu}}=&\int_{-2}^0x \, \sgn{(x+1)(x-b)}dx+1\cdot\sgn{ P'_{b,2}(0)}\\
=&\int_{-2}^{-1}xdx-\int_{-1}^0xdx+\sgn{1-b}=0.
\end{align*}
Then, from Theorem \ref{ThChar-p}, the polynomials $ P_{b,2}$ with $0\leq b < 1$ are all minimal with respect to  \eqref{Example-2} and $\|P_{b,2}\|_{1,\vec{\mu}}=2$.

Furthermore, if $b=1$ the polynomials  $ P_{1,2}(x)=x^2-1$  is minimal and does not satisfy the condition \eqref{sufcond}. Indeed,
\begin{align*}
  \|P_{1,2}\|_{1,\vec{\mu}}= &2 = \|P_{b,2}\|_{1,\vec{\mu}}\quad \text{when }  0\leq b < 1.\\
\langle  P_{1,2} , x\rangle_{1,\vec{\mu}}= & \int_{-2}^0x \, \sgn{x^2-1}dx=-1\neq 0.
\end{align*}
\end{example}

If $1<p<\infty$, from  \cite[Th. 4]{AbelPijeiraLago17}, we know that all minimal  polynomials with respect to \eqref{SobNorm01} (continuous or discrete case) satisfy the condition \eqref{sufcond}.  But as was seen in Example \ref{ExamNoUniqDisc},  this statement is not true when $p=1$. It can even happen that there is no  minimal polynomial  satisfying \eqref{sufcond}.

\begin{example}\label{UniciNoCond}\

Consider the following  discrete Sobolev norm,
\begin{equation}
\|f(x)\|_{1,\vec{\mu}}= \int_{-1}^1 |f(x)|dx+|f'(0)|. \label{discrete_norm_11}
\end{equation}
Then, $P_3(x)=x^3$ is the only 3-th minimal Sobolev polynomial with respect to $\|\cdot\|_{1,\vec{\mu}}$ and does not satisfy the sufficient condition \eqref{sufcond}.

\begin{enumerate}
  \item Note that for every polynomial $Q_n$ we have
$$
\|(-1)^nQ_n(-x)\|_{1,\vec{\mu}}= \int_{-1}^1 |Q_n(-x)|dx+|Q^{\prime}_n(0)|=\|Q_n\|_{1,\vec{\mu}}.
$$
  \item Then, if $S_n$ is a  minimal polynomial  of degree $n$, the monic polynomial $(-1)^nS_n(-x)$ is also extremal. From the convexity of the set of minimal polynomials,
$$P_n(x)=\frac{1}{2} S_n(x)+\frac{(-1)^n}{2}S_n(-x)$$
is an  odd or even polynomial, according to the parity of $n$,  and a monic   minimal polynomial too.
  \item For $n=3$, let  $P_3(x)= x^3+cx$ where $c\in \RR $ a  monic odd polynomial and
\begin{align*}
F(c)=\|x^3+cx\|_{1,\vec{\mu}} = \int_{-1}^1 |x^3+cx|dx + |c| = \left\{\begin{array}{rr}
-2c-\frac{1}{2}, & c \leq-1; \\[.3em]
c^{2}+\frac{1}{2}, &  -1<c<0;\\[.3em]
2c+\frac{1}{2}, & 0 \leq c.
\end{array}\right.
\end{align*}
It is straightforward to see that, the  global minimum of $F$ is attained at $c=0$. Therefore $P_3(x)= x^3$ is a minimal polynomial.
  \item The  polynomial $P_3(x)= x^3$ does not satisfy  \eqref{sufcond}. Indeed,
\begin{align*}
\langle P_3,x\rangle_{1,\vec{\mu}}=\int_{-1}^1 x\, \sgn{x^3}dx=\int_{-1}^1 |x|dx=1\neq 0.
\end{align*}

\item Finally, we will prove the uniqueness. As  $P_3 \in \PP^1_3$ is the only odd  minimal polynomial of degree $3$, and  that  any minimal Sobolev polynomial $S_3 \in  \PP^1_3$ is such that
    $$x^3=\frac{1}{2}S_3(x)-\frac{1}{2}S_3(-x).$$
  Since $\dsty  \|x^3\|_{1,\vec{\mu}}= \frac{1}{2}\|S_3\|_{1,\vec{\mu}}+\frac{1}{2}\|-S_3(-x)\|_{1,\vec{\mu}}$ we get
\begin{align*}
0 \geq \int_{-1}^{1}\left(|x^3|-\frac{1}{2}|S_3(x)|-\frac{1}{2}|S_3(-x)|\right)dx=&|S'_3(0)| \geq 0,
\end{align*}
which implies that   $|x^3|=\frac{1}{2}|S_3(x)|+\frac{1}{2}|S_3(-x)|$ and $|S_3'(0)|=0$. Consequently, $S_3(0)=S_3'(0)=0$ and  $S_3$ takes the form $S_3(x)=x^3+cx^2$, with $c\in\RR$. Since  $c\neq 0$, we arrive at the contradiction
\begin{align*}
\|S_3\|_{1,\vec{\mu}}=\int_{-1}^1|x^3+cx^2|dx=\left\{\begin{array}{ll}
\frac{1}{2}+\frac{1}{6}\,c^4, & |c|< 1; \\[.3em]
\frac{2}{3}|c|, & |c| \geq 1.
\end{array}\right. >\frac{1}{2}=\|x^3\|_{1,\vec{\mu}}.
\end{align*}
So,  $P_3(x)=x^3$ is the only minimal Sobolev polynomial of degree 3.
\end{enumerate}

\end{example}

Note that in this example we have obtained the only monic minimal polynomial of degree $3$ with respect to \eqref{discrete_norm_11}, and it does not satisfy the sufficient condition. This is exclusive to the discrete case. If the vector measure $\vec{\mu}$ is  continuous, the sufficient condition \eqref{sufcond} is also necessary.

\begin{theorem} \label{characterizationP1}
Let  $\vec{\mu}=(\mu_0,\mu_1,\dots,\mu_m)$ be a continuous standard vector measure.
Then, $P_n$ is an $n$-th Sobolev minimal polynomial with respect to $\|\cdot\|_{1,\vec{\mu}}$ if and only if
\begin{equation}\label{neccond}
  \langle P_n , q\rangle_{1,\vec{\mu}}=\sum_{k=0}^m \int_{\Delta_k} q^{(k)}\sgn{P_n^{(k)}}d\mu_k=0, \quad \forall q\in \PP_{n-1}.
\end{equation}
\end{theorem}
\begin{proof}
From Theorem \ref{ThChar-p}, it only remains to prove that the condition \eqref{neccond} is necessary for the extremality. Without loss  of generality, we can assume   that $\deg{P_n}\geq m$, since if $n<m$ we have
$$\|P_n\|_{1,\vec{\mu}} = \sum_{k=0}^n \int_{\Delta_k} \left|P_n^{(k)}\right|\;d\mu_k,\; \text{ and the proof works the same.}$$

Suppose that $P_n\in \PP^1_n$  is a  minimal polynomial with respect to $\|\cdot\|_{1,\vec{\mu}}$  and \eqref{neccond} does not hold. Then there exists  $h\in\PP_{n-1}$ such that  $\langle  P_n , h\rangle_{1,\vec{\mu}}\neq 0$. Multiplying $h$ by a constant we can assume  $\langle P_n , h\rangle_{1,\vec{\mu}}>0$, without loss of generality.

Let $x_{k,1}<x_{k,2}<\dots<x_{k,n_k}$ be the zeros of $P_n^{(k)}$ which lie on $\intersub{\Delta_k}=(a_k,b_k)$. For each $\ell \in \NN$ and $k=0,\dots,m$, denote
\begin{align*}
A_{k,\ell}=\left[a_{k}+\frac{1}{\ell},x_{k,1}-\frac{1}{\ell}\right]\cup \left[x_{k,1}+\frac{1}{\ell},x_{k,2}-\frac{1}{\ell}\right]\cup\cdots \cup\left[x_{k,n_k}+\frac{1}{\ell},b_{k}-\frac{1}{\ell}\right].
\end{align*}
Note, that   $\left\{A_{k,\ell}\right\}_{\ell}$ is a sequence of compact  subsets of  $\intersub{\Delta_k}$, such that  $\dsty  A_k:= \lim_{\ell \to \infty}A_{k,\ell}=\intersub{\Delta_k}\setminus \Lambda_k$, where $\Lambda_k=\{x_{k,1},x_{k,2},\dots,x_{k,n_k}\}$. Let $B_{k,\ell}=\intersub{\Delta_k} \setminus A_{k,\ell}$, so $\dsty \lim_{\ell \to \infty}B_{k,\ell}= \Lambda_k$.

As $\vec{\mu}$ is a vector of  continuous measures,  for every $k=0,1,\dots,m$ we have
\begin{align*}
\lim_{\ell\to \infty}\int_{A_{k,\ell}}h^{(k)}\sgn{P_n^{(k)}}d\mu_k&=\int_{A_{k}}h^{(k)}\sgn{P_n^{(k)}}d\mu_k =\int_{\Delta_k} h^{(k)}\sgn{P_n^{(k)}}d\mu_k,\\
\lim_{\ell\to\infty}\int_{B_{k,\ell}}|h^{(k)}|d\mu_k&=\int_{\Lambda_k}|h^{(k)}|d\mu_k=0.
\end{align*}
Therefore,
\begin{align*}
\lim_{\ell\to \infty}\sum_{k=0}^{m}\int_{A_{k,\ell}}h^{(k)}\sgn{P_n^{(k)}}d\mu_k=&\langle  P_n , h\rangle_{1,\vec{\mu}}>0,\\
\lim_{\ell\to \infty}\sum_{k=0}^{m}\int_{B_{k,\ell}}|h^{(k)}|d\mu_k=&0.
\end{align*}
Hence,  for  $\ell_0 \in \NN$ sufficiently large  $$
\sum_{k=0}^{m}\int_{A_{k,\ell_0}}h^{(k)}\sgn{P_n^{(k)}}d\mu_k>\sum_{k=0}^{m}\int_{B_{k,\ell_0}}|h^{(k)}|d\mu_k.{\cdot}$$
Since every set $A_{k,\ell_0}$, $k=0,1,\dots,m$ is compact and $\Lambda_k \cap  A_{k,\ell_0}=\emptyset$, we get
\begin{align*}
\delta=\min_{k=0,1,\dots,m}\left\{\min_{x\in A_{k,\ell_0}}\{|P_n^{(k)}(x)|\}\right\} >0.
\end{align*}
From the compactness of  $A_{k,\ell_0}$ we also obtain that $$  \delta_h= \max_{k=0,1,\dots,m}\left\{\max_{x\in A_{k,\ell_0}} \left\{|h^{(k)}(x)|\right\}\right\}$$ is   finite and positive. Then we can choose $\lambda>0$ such that $\dsty 0<\lambda<\frac{\delta}{\delta_h}.$

Therefore, for each  $k=0,1,\dots,m$, we have $\dsty |\lambda h^{(k)}(x)|<\delta\leq  |P_n^{(k)}(x)|$ for all $x\in A_{k,\ell_0}$ and
\begin{align*}
\sgn{P_n^{(k)}(x)-\lambda h^{(k)}(x)}=\sgn{P_n^{(k)}(x)}, \quad  \text{for all }   x\in A_{k,\ell_0}.
\end{align*}
Finally,
\begin{align*}
\|P_n-\lambda h\|_{1,\vec{\mu}}=&\sum_{k=0}^{m}\int_{\Delta_k} |P_n^{(k)}-\lambda h^{(k)}|d\mu_{k} \\
=&\sum_{k=0}^{m}\left(\int_{B_{k,\ell_0}}\!\!\!\!\!|P_n^{(k)}-\lambda h^{(k)}|d\mu_k+\int_{A_{k,\ell_0}}\!\!\!\!\!|P_n^{(k)}-\lambda h^{(k)}|d\mu_k\right)
\end{align*}
\begin{align*}
=&\sum_{k=0}^{m}\left(\int_{B_{k,\ell_0}}\!\!\!\!\!|P_n^{(k)}-\lambda h^{(k)}|d\mu_k + \int_{A_{k,\ell_0}}\!\!\!\!\!\sgn{P_n^{(k)}-\lambda h^{(k)}}\left(P_n^{(k)}-\lambda h^{(k)}\right)d\mu_k\right) \\
=&\sum_{k=0}^{m}\left(\int_{B_{k,\ell_0}}\!\!\!\!\!|P_n^{(k)}-\lambda h^{(k)}|d\mu_k  + \int_{A_{k,\ell_0}}\!\!\!\!\!\sgn{P_n^{(k)}}\!\!\left(P_n^{(k)}-\lambda h^{(k)}\right)d\mu_k\right) \\
\leq & \sum_{k=0}^{m}\left(\int_{B_{k,\ell_0}}\!\!\!\!\!|P_n^{(k)}|d\mu_k+\lambda\int_{B_{k,\ell_0}}\!\!\!\!\!|h^{(k)}|d\mu_k \right. \\ & + \left. \int_{A_{k,\ell_0}}\!\!\!\!\!|P_n^{(k)}|d\mu_k-\lambda\int_{A_{k,\ell_0}}\!\!\!\!\!\sgn{P_n^{(k)}}h^{(k)}d\mu_k\right) \\
=&\sum_{k=0}^{m}\int_{\Delta_k} |P_n^{(k)}|d\mu_{k} +\lambda\left(\sum_{k=0}^{m}\int_{B_{k,\ell_0}}\!\!\!\!\!|h^{(k)}|d\mu_k-\sum_{k=0}^{m}\int_{A_{k,\ell_0}}\!\!\!\!\!\sgn{P_n^{(k)}}h^{(k)}d\mu_k\right) \\
< & \|P_n\|_{1,\vec{\mu}},
\end{align*}
which is a contradiction with the extremality of $P_n$.
\qed \end{proof}

\section{Lacunary and non-lacunary discrete  Sobolev norms} \label{Sec3}

Most of the formulas given here are known to the specialist, although precise references may be hard to find in the literature.  Therefore, we include this  section with full proofs for completeness, except when an exact reference is available.

Consider a finite positive Borel measure $\mu$, being  $\supp{\mu}$ a subset of the real line with infinitely many points such that  $\PP\subset L^{1}(\mu)$. In the remainder,  we assume  that $N \in \ZZp$, $\Omega=\{ c_1,c_2,\dots,c_N\}\subset \RR $, $\{m_0,m_1,\dots,m_N\} \subset \ZZp$ and $m=\max\{m_0,m_1,\dots,m_N\}$. Let  $\vec{\mu}=(\mu_0,\mu_1,\dots,\mu_m)$ be the  standard vector measure. For each $ 1 \leq p<\infty$, let us consider the general discrete Sobolev norm
\begin{align}
\|f\|_{p,\vec{\mu}}= & \left(\sum_{k=0}^m \int_{\Delta_k}  \left|f^{(k)}\right|^p d\mu_k\right)^{1/p} =\left(\int_{\Delta} \left|f\right|^{p}\;d\mu+ \sum_{j=1}^N \sum_{k=0}^{m_j} A_{j,k}\left|f^{(k)}(c_j)\right|^p\right)^{1/p},\label{discreteSnorm}
\end{align}
 where $\Delta$ is  the convex hull of the support of the measure $\mu$. Notice that, unlike \eqref{SobNorm01}, the representation \eqref{discreteSnorm} of $\|\cdot\|_{p, \vec{\mu}}$ is not unique, but depends on how many Dirac measures, of the discrete part of $\mu_0$, are included in the measure $\mu$.   In general, the representation \eqref{discreteSnorm} is unique once the measure $\mu$ is fixed, so this dependence will be omitted for brevity.

If   there exists a constant $M$ such that \begin{equation}\label{MUltOperato}
\|xq\|_{p,\vec{\mu}} \leq M  \|q\|_{p,\vec{\mu}}, \quad \text{for all } q \in \PP,
\end{equation}  we say that the multiplication operator is bounded on  $\PP$ with respect to $\|\cdot\|_{p,\vec{\mu}}$. The close  relation between \eqref{MUltOperato} and the uniform boundedness of the set of zeros of sequences of minimal polynomials was established   in \cite{HPLago99}. Since then, several studies have been published on this subject.

\begin{proposition} \label{PropoNLacu-02}Assume that the discrete Sobolev norm \eqref{discreteSnorm}   is non-lacunary and $\Delta$ is bounded,
then for each $q  \in \PP$ we have
\begin{equation*}
 \|xq\|_{p,\vec{\mu}} \leq M\|q\|_{p,\vec{\mu}},
\end{equation*}
\begin{align*}
\text{where }  \; M= &\max\left\{M_1,2^{p-1}(M_1+mM_2)\right\}^{1/p}, \;  M_1=\sup_{x\in K} |x|^p, \;   K = \Delta \cup \{c_1,\ldots,c_m\} , \\
   M_2= & \max \left\{\frac{A_{j,k+1}}{A_{j,k}}:1\leq j \leq N \text{ and } 0 \leq k \leq m_j-1 \right\}.
\end{align*}
\end{proposition}

\begin{proof}
Notice that $(x q)^{(k)}=x q^{(k)}+k q^{(k-1)}, \quad k\in\NN.$ Therefore
\begin{align*}
\Psi:=& \sum_{j=1}^N \sum_{k=0}^{m_j} A_{j,k}\left|c_j q^{(k)}(c_j)+k q^{(k-1)}(c_j)\right|^p\\
\leq & 2^{p-1}\left( \sum_{j=1}^N \sum_{k=0}^{m_j} A_{j,k}\left| c_j q^{(k)}(c_j)\right|^p + \sum_{j=1}^N \sum_{k=1}^{m_j} A_{j,k} \left| k q^{(k-1)} (c_j)\right|^p\right)
\\
\leq & 2^{p-1}\left(  M_1 \sum_{j=1}^N \sum_{k=0}^{m_j} A_{j,k}\left|q^{(k)}(c_j)\right|^p     +m \sum_{j=1}^N \sum_{k=1}^{m_j} A_{j,k} \left|q^{(k-1)}(c_j)\right|^p\right)
\\
= &  2^{p-1}\left(  M_1 \sum_{j=1}^N \sum_{k=0}^{m_j} A_{j,k}\left|q^{(k)}(c_j)\right|^p   + m\sum_{j=1}^N \sum_{k=0}^{m_j-1} A_{j,k+1} \left|q^{(k)}(c_j) \right|^p\right)\\
\leq & 2^{p-1}\left(  M_1\sum_{j=1}^N \sum_{k=0}^{m_j} A_{j,k}\left| q^{(k)}(c_j)\right|^p    +mM_2\sum_{j=1}^N \sum_{k=0}^{m_j-1} A_{j,k} \left|q^{(k)}(c_j)\right|^p\right) \\
\leq & 2^{p-1}\left(  (M_1+mM_2)\sum_{j=1}^N\sum_{k=0}^{m_j} A_{j,k}\left| q^{(k)}(c_j)\right|^p\right).\\
\|xq\|^{p}_{p,\vec{\mu}}= &  \int_{ \Delta}|xq|^p d \mu + \Psi \\ \leq  & M_1\int_{ \Delta}|q|^p d\mu  + 2^{p-1}\left(  (M_1+mM_2)\sum_{j=1}^N\sum_{k=0}^{m_j} A_{j,k}\left| q^{(k)}(c_j)\right|^p\right)  \leq     M^{p}\|q\|^p_{ p,\vec{\mu}}.
\end{align*}
\qed \end{proof}

If  $\|\cdot\|_{p,\vec{\mu}}$ is a lacunary Sobolev norm defined as in \eqref{discreteSnorm}, we  define the \emph{associated non-lacunary norm} as  $\|\cdot\|_{p,\vec{\mu}^*}$
\begin{equation}\label{discreteSnorm-nL}
\|f\|_{p,\vec{\mu}^*}= \left(\int_{\Delta} \left|f\right|^{p}\;d\mu+ \sum_{j=1}^N \sum_{k=0}^{ m_j} A_{j,k}^*\left|f^{(k)}(c_j)\right|^p\right)^{1/p}.
\end{equation}
where $\dsty A_{j,k}^*=\left\{
                     \begin{array}{ll}
                       A_{j,k}, & \hbox{if } A_{j,k}>0 \text{ or }m_j<k \leq m; \\
                       1, &  \hbox{in other case}.
                     \end{array}
                   \right.$

\begin{proposition} \label{PropoNLacu-03} Let $\|\cdot\|_{p,\vec{\mu}}$ be a lacunary Sobolev norm defined as in \eqref{discreteSnorm}, with  $\Delta$  bounded. Then, there exists a constant $M$ such that
$\|x q\|_{p,\vec{\mu}} \leq M \|q\|_{p,\vec{\mu}}$ for all $q \in \PP$ if and only if the lacunary norm \eqref{discreteSnorm} and the associated non-lacunary norm \eqref{discreteSnorm-nL} are  equivalents (i.e. $\|\cdot\|_{p,\vec{\mu}}$ is essentially non-lacunary).
\end{proposition}

\begin{proof}  Assume that  a lacunary norm defined as in \eqref{discreteSnorm} is equivalent to its associated non-lacunary norm \eqref{discreteSnorm-nL}.  From   Proposition \ref{PropoNLacu-02}, it is  straightforward that there exists a constant $M$ such that $\|x q\|_{p,\vec{\mu}} \leq M \|q\|_{p,\vec{\mu}} $.

Now, suppose that the multiplication operator is bounded on  $\PP$ with respect to  the lacunary norm $\|\cdot\|_{p,\vec{\mu}}$, then there exist $M>0$ : $\|xq\|_{p,\vec{\mu}}\leq \|q\|_{p,\vec{\mu}}$, $q\in \PP$. From \eqref{discreteSnorm-nL}, obviously  $\|q\|_{p,\vec{\mu}} \leq  \|q\|_{p,\vec{\mu}^*}$.  Furthermore, from definition
\begin{align*}
  \|q\|_{p,\vec{\mu}^*}= &  \left(\|q\|_{p,\vec{\mu}}^p+ \sum_{j=1}^N \sum_{k\in I_j} \left|q^{(k)}(c_j)\right|^p\right)^{1/p}
 \leq \|q\|_{p,\vec{\mu}}+ \left(\sum_{j=1}^N \sum_{k\in I_j} \left|q^{(k)}(c_j)\right|^p\right)^{1/p}\ ,
\end{align*}
where $I_j=\{k: A_{j,k}=0 \text{ and } 0\leq k < m_j\}$. Therefore, the remainder of the proof is devoted to find a constant $K^*$ such that
\begin{equation}\label{Inequality-03}
  \left( \sum_{j=1}^N \sum_{k\in I_j} \left|q^{(k)}(c_j)\right|^p\right)^{1/p}\leq K^{*} \;   \|q\|_{p,\vec{\mu}} \quad q\in\PP\,.
\end{equation}
To achieve this purpose,  it is sufficient to prove that for  every  $j$ and $0\leq k <  m_j$ there exists a constant $K_{j,k}>0$ satisfying
\begin{equation}\label{Inequality-01}
\left|q^{(k)}(c_j) \right| \leq     K_{j,k} \|q\|_{p,\vec{\mu}} \quad q\in\PP\, .
\end{equation}
In this case,  taking $ \dsty K^{*} = \left( \sum_{j=1}^N \sum_{k\in I_j} K_{j,k}^p\right)^{1/p}$, we get   \eqref{Inequality-03}.

To prove the  inequality \eqref{Inequality-01}, note that
\begin{align} \nonumber
\left|(k+1) q^{(k)}(c_j)\right|- \left|c_j q^{(k+1)}(c_j)\right| \leq & \left|(k+1) q^{(k)}(c_j)+c_j q^{(k+1)}(c_j)\right|= \left|(xq)^{(k+1)}(c_j)\right|,\\
\nonumber
\left| q^{(k)}(c_j)\right| \leq \left|(k+1) q^{(k)}(c_j)\right| \leq &\left|(xq)^{(k+1)}(c_j)\right| +  \left|c_j q^{(k+1)}(c_j)\right|
\\ \leq & \left|(xq)^{(k+1)}(c_j)\right| + |c^*| \left| q^{(k+1)}(c_j)\right|,\label{Inequality-02}
\end{align}
where $\dsty c^*= \max_{1\leq j \leq N}|c_j|$.  If $m_j-k=1$, and $q\in\PP$
\begin{align*}
 \left| q^{(m_j-1)}(c_j)\right| \leq &\frac{1}{A_{j,m_j}}  \left|A_{j,m_j} (xq)^{(m_j)}(c_j)\right| + \frac{|c^*|}{A_{j,m_j}}  \left|A_{j,m_j}q^{(m_j)}(c_j)\right|. \\
 \leq  & \frac{1}{A_{j,m_j}}  \|xq\|_{p,\vec{\mu}} + \frac{|c^*|}{A_{j,m_j}} \|q\|_{p,\vec{\mu}} \leq  K_{j,m_j-1}\; \|q\|_{p,\vec{\mu}}.
\end{align*}
where $\dsty K_{j,m_j-1}= \frac{M+|c^*|}{A_{j,m_j}} \neq 0$ and we get \eqref{Inequality-01} for $k= m_j-1$.

We now proceed by induction.
\begin{enumerate}
  \item {$[m_j-k=\ell]$} Assume that \eqref{Inequality-01} holds for $k=m_j-\ell$, i.e. there exists a constant $ K_{j,m_j-\ell}\neq 0$ such that
  $$  \left| q^{(m_j-\ell)}(c_j)\right| \leq  K_{j,m_j-\ell}\; \|q\|_{p,\vec{\mu}}.  $$

  \item {$[m_j-k=\ell+1]$} If   $k=m_j-\ell-1$, from \eqref{Inequality-02} and the induction hypothesis
  \begin{align*}
    \left| q^{(m_j-\ell-1)}(c_j)\right| \leq & \left|(xq)^{(m_j-\ell)}(c_j)\right| + |c^*| \left| q^{(m_j-\ell)}(c_j)\right| \\ \leq & K_{j,m_j-\ell}\; \|xq\|_{p,\vec{\mu}}+K_{j,m_j-\ell}\;|c^*| \, \|q\|_{p,\vec{\mu}}   \leq    K_{j,m_j-\ell-1}\;\|q\|_{p,\vec{\mu}},
  \end{align*}
\end{enumerate}
  where $K_{j,m_j-\ell-1}= (M+|c^*|) K_{j,m_j-\ell}$. \qed \end{proof}

\begin{theorem}\label{ThEssNonL} If \eqref{discreteSnorm} is essentially non-lacunary,   then the set of zeros of a minimal polynomial sequence is  uniformly bounded.
\end{theorem}

\begin{proof}  Let \eqref{discreteSnorm} be  an essentially non-lacunary Sobolev norm and \eqref{discreteSnorm-nL} its  associated non-lacunary Sobolev norm. From Proposition \ref{PropoNLacu-03}, there exist   constants $C_1,C_2>0$ such that $C_1\,\|q\|_{p,\vec{\mu}^*}\leq  \|q\|_{p,\vec{\mu}} \leq  C_2\,\|q\|_{p,\vec{\mu}^*}$ for all $q \in \PP$.  Moreover, from Proposition \ref{PropoNLacu-02}, there exists another constant $C_3>0$ such that $ \|z\,q\|_{p, \vec{\mu}^*} \leq  C_3\,\|q\|_{p,\vec{\mu}^*}$ .

If   $P_n$ is a  minimal polynomial of degree $n$ and   $P_n(z_0)=0$,  there exists a monic polynomial $q$ of degree $n-1$ such that $P_n(z)= (z-z_0)q(z)$. As $P_n$ is minimal $$\left|z_{0}\right|\|q\|_{p,\vec{\mu}}-\|z q\|_{p,\vec{\mu}} \leq\left\|z_{0} q-z q\right\|_{p,\vec{\mu}}=\left\|P_{n}\right\|_{p,\vec{\mu}} \leq\|z q\|_{p,\vec{\mu}}.$$
Then,
$$ \left|z_{0}\right|C_1\|q\|_{p,\vec{\mu}^*}  \leq \left|z_{0}\right|\|q\|_{p,\vec{\mu}} \leq 2\|z q\|_{p,\vec{\mu}} \leq 2C_2\|z q\|_{p,\vec{\mu}^*}\leq 2C_2C_3 \| q\|_{p,\vec{\mu}^*},$$
which completes the proof.
\qed \end{proof}

\subsection{Asymptotic distribution of zeros}

To state the result on the zero distribution of minimal polynomials with respect to an  essentially non-lacunary  norm, we need to  introduce some concepts and notations.

\begin{itemize}
  \item For any polynomial $q$ of exact degree $n$, we denote $\dsty \vartheta(q)= \dfrac{1}{n}\sum_{j=1}^n\delta_{z_j},$  where $z_1,\dots,z_n$ are the zeros of $q$ repeated according to their multiplicity. This is the so called normalized counting measure associated with $q$.
  \item If  $\Delta=\supp\mu$  is regular (a compact subset of the complex plane is said to be regular if the unbounded connected component of its complement is regular with respect to the Dirichlet problem), the measure $\mu\in\mathbf{Reg}$ if and only if
\begin{equation}\label{RegCond}
\lim_{n\to\infty} \left({\frac{\|q_n\|_{\Delta}}{\|q_n\|_{ p,\mu }}}\right)^{1/n} = 1,
\end{equation}
for every sequence of polynomials $\{q_n\}$, $\deg{q_n}\leq n$, $q_n\not\equiv 0$  (cf. \cite[Th 3.4.3]{Stahl92}), where  $\|\cdot\|_{\mathcal{A}}$ denotes the supremum norm on $\mathcal{A}\subset \CC$.
  \item   Given a compact set $\mathcal{A} \subset \CC$,  $\cp{\mathcal{A}}$ denotes the logarithmic capacity of $\mathcal{A}$,  $\omega_{\mathcal{A}}$ the equilibrium measure on $\mathcal{A}$ and $\displaystyle G_{\mathcal{A} }(z;\infty)$ the corresponding Green's function with singularity at infinity (cf.\cite{Ran95,Stahl92}).
  \item Let $T_n$ be the $n$-th  monic minimal  polynomial with respect to $\|\cdot\|_{\Delta}$, i.e.  the $n$-th Chebyshev polynomial with respecto to $\Delta$. It is known that
      \begin{equation}\label{ChebyConstant}
\lim_{n\to\infty} \|T_n\|^{1/n}_{\Delta}= \cp{\Delta}. \quad \text{\cite[Cor. 5.5.5]{Ran95}}
      \end{equation}

\end{itemize}

To determine the asymptotic distribution of zeros of sequences of minimal  polynomials   in this section, we need the following lemma.

\begin{lemma}{\cite[Lemma 3]{HPLago99}}
Let $E$ be a compact regular subset of the complex plane and $\displaystyle\lbrace q_n\rbrace$ a sequence of polynomials such that $\deg{q_n}\leq n$ and $q_n\not\equiv 0$. Then, for all $k\in \mathbb{Z}_{+}$,
\begin{equation}\label{lemma3}
\uplim_{n\to\infty}\sqrt[n]{ \frac{\|q_n^{(k)}\|_E}{\|q_n\|_E}}\leq 1.
\end{equation}
\end{lemma}

The following theorem is the main result of  this section and is valid for discrete Sobolev norms, whether lacunary or not.
For   $p=2$,   the theorem was proved in \cite[Th. 5]{HPLago99},  and for continuous Sobolev norms in \cite[Th. 2]{GlIpHp05}. The scheme of the proof is quite similar to the previous ones.

\begin{theorem}\label{AsymZerosDis}
Consider a  discrete Sobolev $p$-norm \eqref{discreteSnorm}, such that $\mu\in\mathbf{Reg}$ and $\Delta$ is a bounded real interval.
If  $\{P_n\}$ is the sequence of monic minimal polynomials with respect to \eqref{discreteSnorm}, then  for all $j\geq  0$

\begin{align}\label{eq8.1}
\lim_{n\to\infty}\|P_n^{(j)}\|_{\Delta}^{1/n}= &\cp{\Delta}, \quad \text{and} \\
\label{eq9.1}
\wlim_{n\to\infty}\vartheta\left(P_n^{(j)}\right) =&\omega_{\Delta}, \quad \text{in the weak topology of measures.}
\end{align}

\end{theorem}

\begin{proof}
Firstly, the compact set $\Delta$ has empty interior and connected complement and under these conditions (see \cite[Th. 2.1]{Blatt88}) we have that \eqref{eq8.1} implies \eqref{eq9.1}.

Let $T_n$ be the $n$-th monic minimal  polynomial with respect to $\|\cdot\|_{\Delta}$, i.e.  the $n$-th Chebyshev polynomial with respecto to $\Delta$. From \eqref{ChebyConstant}, it is straightforward to see that for all sequence $\{ Q_n \}_{n\in\mathbb{Z}_+}$ of monic polynomials $Q_n$ of degree $n$

\begin{equation}\label{eq8.1-1}
\lowlim_{n\to\infty} \|Q_n^{(j)}\|_{\Delta}^{1/n} \geq \lowlim_{n\to\infty} \|T_{n-j}\|_{\Delta}^{1/n} = \cp{\Delta}.
\end{equation}

 If $ \rho(z)= \prod_{j=1}^N (z-c_j)^{m_j+1}$   and  $n\geq  \mathbf{d}:=N+\sum_{j=1}^{N}\!m_j$, we get

\begin{align*}
\|P_n\|^{ p}_{p,\mu} & \leq  \|P_n\|^{p}_{p,\vec{\mu}}\leq  \| \rho\,T_{n-\mathbf{d}}\|^{p}_{p,\vec{\mu}} =\int_{\Delta}\left| \rho\,T_{n-\mathbf{d}}\right|^p d\mu\leq \mu\left(\Delta\right)\|\rho\|^p_{\Delta}\|T_{n-\mathbf{d}}\|^p_{\Delta}.
\end{align*}

From \eqref{RegCond}-\eqref{ChebyConstant},     $\dsty \uplim_{n\to\infty}\|P_n\|^{1/n}_{\Delta}\leq \cp{\Delta}$. Therefore, as $\Delta$ is a   compact regular set, from  \eqref{lemma3}  we have  for  every  $j \geq  0$

\begin{equation}\label{eq8.1-2}
\uplim_{n\to\infty}\|P_n^{(j)}\|_{\Delta}^{1/n}\leq\cp{\Delta}.
\end{equation}

Finally, from \eqref{eq8.1-1}-\eqref{eq8.1-2} we get  \eqref{eq8.1}. \qed

\end{proof}

If the norm \eqref{discreteSnorm} is essentially non-lacunary, from  Theorem \ref{ThEssNonL},   we know that  there exists a constant $M$ such that

$$\{z \in \CC: P_n(z)=0 \; \text{for some } n \in \ZZp\} \subset D_M=\{z \in \CC: \; |z|\leq M\},$$
where $\{P_n\}$ is a sequence of minimal polynomials with respect to \eqref{discreteSnorm} ($\dgr{P_n}=n$). Under this consideration we have the following asymptotic results.

\begin{corollary}
Assume that $\{P_n\}$ is the sequence of minimal polynomials with respect to an essentially non-lacunary  norm \eqref{discreteSnorm}, where   $\Delta$ is regular and $\mu\in\mathbf{Reg}$.  Then, for all $j\in\mathbf{Z}_+$
\begin{enumerate}
  \item $\dsty \uplim_{n\to\infty}\left|P_n^{(j)}(z)\right|^{1/n}=\cp{\Delta} e^{G_{\Delta}(z;\infty)}, $ for every $z\in\CC$ except for a set of capacity zero,
  \item $\dsty \lim_{n\to\infty}\left|P_n^{(j)}(z)\right|^{1/n}=\cp{\Delta}e^{G_{\Delta}(z;\infty)}$,  uniformly on compact subsets of $\Omega= \CC \setminus D_M$.
  \item $\dsty \lim_{n\to\infty}\frac{P_n^{(j+1)}(z)}{nP_n^{(j)}(z)}=\int_{\Delta}\frac{d\omega_{\Delta}(x)}{z-x}$,  uniformly on compact subsets of $\Omega$.
\end{enumerate}
\end{corollary}

\begin{proof} From Proposition \ref{PropoNLacu-03}, it is sufficient to prove the corollary for non-lacunary norms. As it was commented for the case  $p=2$ in the last paragraph of \cite{HPLago99}, the proof here follows \cite[Th. 6]{GlIpHp05} point by point to get  the desired result.
\qed \end{proof}

\section{Sequentially-ordered discrete Sobolev norm}\label{Sec4}

If the discrete Sobolev norm \eqref{discreteSnorm} is non-lacunary, it is easy to prove that the $n$-th minimal Sobolev polynomial has all its the zeros located on $\Delta$,  except a number of them equal to the amount of  non-zero values $A_{j,k}$ in the discrete part of \eqref{discreteSnorm};  see Proposition \ref{ZeroLoc-01}. In this section, we extend this result to lacunary Sobolev norms  when the discrete part of \eqref{discreteSnorm} satisfies certain order condition.

Fix $1<p<\infty$ and a standard vector measure $\vec{\mu}$ such that $\|\cdot\|_{p,\vec{\mu}}$ is a discrete Sobolev norm defined by \eqref{discreteSnorm} and satisfying $c_j\notin \inter{\Delta}=(a,b)$ for $j=1,2,\dots,N$.  As in the previous section,  consider the polynomial
	$$\rho(x)= \prod_{c_j\leq a}\!\!\left(x-c_j\right)^{m_{j}+1} \!\prod_{c_j\geq b}\!\left(c_j-x\right)^{m_{j}+1}$$
of degree $ \mathbf{d}=N+\sum_{j=1}^{N}\!m_j$ and positive on $(a,b)$.
If $n > \mathbf{d}$ and $P_n$ is the $n$-th minimal polynomial with respect to \eqref{discreteSnorm},  from Theorem  \ref{TheoremPProduct}
\begin{equation}\label{quasi-orthogonal}
	\int_{a}^b q \, \sgn{P_n}|P_n|^{p-1}\rho \,d\mu=  \langle  P_n,q\rho\rangle_{p,\vec{\mu}}= 0 ,
\end{equation}
for every $q\in \PP_{n- \mathbf{d}-1}$. Hence, the polynomial $P_n$ hast at least $n- \mathbf{d}$ changes of sign on $\inter{\Delta}$, otherwise  \eqref{quasi-orthogonal} lead us to a contradiction with
$$	\int_{a}^b q \, \sgn{P_n}|P_n|^{p-1}\rho d\mu>0,
$$
where $q$ is the polynomial having a simple zero on each change of sign of $P_n$ on $(a,b)$.  So, we have proved the following proposition, which is the extension of \cite[Proposition 2.1]{AbelPijeiraIgna20} to the minimal case, $1<p<\infty$.
\begin{proposition}\label{ZeroLoc-01}
	Let $P_n$ be the $n$-th  Sobolev minimal polynomial with respect to \eqref{discreteSnorm} ($1<p<\infty$),  which satisfies $c_j\notin \inter{\Delta}$ for $j=1,2,\dots,N$, and  $n>  \mathbf{d}$, then $P_n$ has at least $(n- \mathbf{d})$ changes of sign on $\inter{\Delta}$.
\end{proposition}

Proposition \ref{ZeroLoc-01} can also be seen as a generalization of the zero location theorem for standard orthogonal polynomials ($p=2$ and $m=0$). However,  a result proved by M. G. Bruin already in 1993, see \cite[Th. 4.1]{Bruin93}, seems to suggest that the number of zeros of $P_n$ in $\inter{\Delta}$ does not depend only on the higher order derivatives $m_j$ of each point $c_j$, but on the number of terms in the discrete part of \eqref{discreteSnorm}
	$$ \mathbf{d}^*:=\left|\{A_{j,k}>0: j=1,2,\dots,N, \ k = 0, 1,\dots , m_j\}\right|,$$
where $|A|$ denotes the cardinality of a set $A$.

 This assumption became even stronger when the relative asymptotic of discrete Sobolev orthogonal polynomials \cite[Theorem 4]{LagMarWal95} was found. Finally, in \cite{AlLagRez96}, the authors  proved it for the case when  \eqref{discreteSnorm} has only one mass point ($N=1$).

\begin{theorem}[{\cite[Th. 2.2]{AlLagRez96}}]\label{ThZeroLocDiscreteLago}
Let $\mu$ be a standard measure such that $c\in\RR\setminus \inter{\Delta}$. If $P_n$ denotes the $n$-th Sobolev minimal polynomial with respect to
\begin{align*}
	\|f\|_{2,\vec{\mu}} =\; \left(\int_{\Delta} |f|^2d\mu+\sum_{k=0}^{m}A_{k} |f^{(k)}(c)|^2\right)^{1/2}.
\end{align*}
Then $P_n$ has at least $n- \mathbf{d}^*$ changes of sign in $\inter{\Delta}$.
\end{theorem}

 The next examples show that this theorem is not longer true if we consider arbitrary mass point configurations with more than one point (i.e.  $N \geq2$ in \eqref{discreteSnorm}), at least not for every value of $n$.

\begin{example}[bounded case]\label{ExampleNoSeqOrdAcot}
	Set $$ \|f\|_{2,\vec{\mu}} = \left( \int_{-1}^{1} |f|^2 dx+8|f^{\prime}(4)|^2+6|f^{\prime\prime}(2)|^2\right)^{1/2},$$		then
	\begin{align*}
		P_4(x)= k_4\!\left(x^4-\frac{2595}{803}x^3-\frac{5232}{539}x^2-\frac{837735}{39347}x+\frac{8181}{2695}\right),
	\end{align*}
	whose zeros are approximately $\xi_{1}\approx 0.13 $, $\xi_{2}\approx-5.62$, $\xi_{3}\approx -1.26+1.56i$ and $\xi_{4}\approx -1.26-1.56i$.
\end{example}

\begin{example}[unbounded case]\label{ExampleNoSeqOrdNoAcot}
	Set $$ \|f\|_{2,\vec{\mu}}= \left( \int_{0}^{\infty} |f(x)|^2e^{-x} dx+3|f^{\prime}(-4)|^2 +8|f^{\prime\prime}(0)|^{2}\right)^{1/2},$$  then
	\begin{align*}
		P_4(x)= k_4\!\left(x^4-\frac{128}{97}x^3-\frac{2536}{97}x^2+\frac{8800}{97}x-\frac{5288}{97}\right),
	\end{align*}
	whose zeros are approximately $\xi_{1}\approx 0.78 $, $\xi_{2}\approx-5.93$, $\xi_{3}\approx 3.24+1.16i$ and $\xi_{4}\approx 3.24-1.16i$.
\end{example}

Note that, in both cases, three zeros of $P_4$ are out of $\inter{\Delta}$ and two of them  are non-real.

The first result treating the case $N\geq 2$ in a general way is \cite[Theorem 1]{AbelPijeiraIgna20}. Here, the authors give a result similar to Theorem \ref{ThZeroLocDiscreteLago} for $N\geq 2$ in the case  $p=2$ and the discrete part of  \eqref{discreteSnorm}  satisfies certain order condition. The condition was called by the authors   the sequentially order condition. Although the condition was enough for the purposes of the paper,  it does not include the case of Theorem \ref{ThZeroLocDiscreteLago}, when there is more than one order derivative at the same mass point $c_j$. Following the same technique,    we expand this condition a little bit more,   in such a way that the case of Theorem \ref{ThZeroLocDiscreteLago} is included. We will remain calling it the sequentially order condition or we  will simply say that the discrete Sobolev norm is sequentially ordered. The result is also generalized for the minimal case $1<p<\infty$.

\begin{definition}[\textbf{Sequentially-ordered Sobolev norm}]
	We say that a discrete Sobolev norm $\|\cdot\|_{p,\vec{\mu}}$ defined by \eqref{discreteSnorm},  is \emph{sequentially ordered} if the conditions
	\begin{equation*}
\Delta_{k}\cap \intch{\cup_{i=0}^{k-1}\Delta_i} =\emptyset, \quad \quad k=1,2,\dots, m,\quad \text{hold.}
	\end{equation*}		

\end{definition}
We recall that $\Delta_k:=\ch{\supp{\mu_k}}$, so in the discrete case they can be rewritten as
$$\Delta_k=\funD{\ch{\Delta\cup\{c_j: A_{j,0}>0\}}}{k=0}{\ch{\{c_j: A_{j,k}>0\}}}{1\leq k\leq m.}$$

\begin{example} The following Sobolev discrete norms are sequentially ordered for any $p\in[1,\infty)$ and a standard measure $\mu$
\begin{align*}
\|f\|_{p,\vec{\mu}}=&\left(\int_{-1}^{1}|f|^pd\mu+4|f^{\prime}(-1)|^p+|f^{\prime}(-3)|^{p}+3|f^{\prime\prime}(2)|^{p}+5|f^{(5)}(-3)|^{p}\right)^{1/p}.\\
	\|f\|_{p,\vec{\mu}}=&\left(\int_{-1}^{1}|f|^pd\mu+\sum_{k=0}^{\ell_1}A_{1,k}|f^{(k)}(-1)|^p+\sum_{k=0}^{\ell_2}A_{2,k} |f^{(k)}(1)|^p \right)^{1/p}.
\end{align*}
where $A_{1,k}A_{2,k}=0$ for $k=0,1,\dots,\min\{\ell_1,\ell_2\}$.
\end{example}

\begin{theorem}\label{ThZeroLocDiscrete}
Let $\vec{\mu}$ be a standard vector measure and $1<p<\infty$. If $\|\cdot\|_{p,\vec{\mu}}$ is a  sequentially-ordered Sobolev norm written as  \eqref{discreteSnorm}, where $\mu$ is taken in such a way $c_j\notin \inter{\Delta}$,  then  $P_n$ has at least $n-\mathbf{d}^*$ changes of sign on $\inter{\Delta}$.
\end{theorem}

It is worth noting that, although the theorem is enunciated depending on which representation \eqref{discreteSnorm} of the Sobolev norm is considered, the definition of sequentially ordered Sobolev norm   is independent of this representation. If what we are after  is to   locate  the largest possible  number of zeros, we should calculate $\mathbf{d}^*$ in the theorem considering the representation \eqref{SobNorm01}, rather than \eqref{discreteSnorm}. However,   in this case we would have the zeros located in the bigger set $\Delta_0\supset \Delta$. Because of the assumption $c_j\notin\inter{\Delta}$, this inclusion is strict except for the trivial case of \eqref{SobNorm01} and \eqref{discreteSnorm} agree ($\mu\equiv \mu_0$).

Notice that both Examples \ref{ExampleNoSeqOrdAcot} and \ref{ExampleNoSeqOrdNoAcot} are not sequentially-ordered. So, this order restriction in the discrete part seems to be optimal to have the most number of zeros simple and located on $\inter{\Delta}$, at least for every value of $n$.

\subsection{Proof of Theorem \ref{ThZeroLocDiscrete}}

Given a polynomial $Q$ with real coefficients and a real set $A$, we introduce the following notations:
\begin{itemize}
			\item $\Ncero{Q}{A}$ denotes the number of values on $A$ where the polynomial $Q$ vanishes, (i.e. zeros of $Q$ on $A$ without counting multiplicities).
			\item $\Nzero{Q}{A}$ denotes the total number of zeros (counting multiplicities) of $Q$ on $A$.
\end{itemize}

The next lemma is an extension of \cite[Lem. 2.1]{LagIpHp01} and \cite[Lem. 3.1]{AbelPijeiraIgna20}.
\begin{lemma}
	Let $\{I_{k}\}_{k=0}^{m}$ be a set of intervals on the real line with $m\in\ZZ_+$ and let $Q$ be a polynomial with real coefficients of degree $\geq m$. If
	\begin{align}\label{CondSeqOrdSet}
	I_{k}\cap \intch{\cup_{i=0}^{k-1}I_i}=\emptyset, \quad \quad k=1,2,\dots, m,
	\end{align}
	then
	\begin{align}\nonumber
		\Nzero{Q}{J}+\Ncero{Q}{I_0\setminus J}+\sum_{i=1}^m \Ncero{Q^{(i)}}{I_i}\leq & \Nzero{Q^{(m)}}{J}\\ \label{IneqRollGeneral}
		 & + \Ncero{Q^{(m)}}{\ch{\cup_{i=0}^{m}I_i}\setminus J}+m,
	\end{align}
	for every closed subinterval $J$ of $\intersub{I_0}$ (both empty set and unitary sets are assumed to be intervals).
\end{lemma}
\begin{proof}
	First, we are going to point out the following consequence of Rolle's Theorem. If $I$ is a real interval and $J$ is a closed subinterval of $\inter{I}$, then
	\begin{align}\label{IneqRoll}
		\Nzero{Q}{J}+\Ncero{Q}{I\setminus J}\leq \Nzero{Q'}{J}+ \Ncero{Q'}{\inter{I}\setminus J}+1.
	\end{align}	
	For $m=0$ \eqref{IneqRollGeneral} trivially holds. We now proceed by induction on $m$. Suppose that we have  $m+1$ intervals  $\{I_i\}_{i=0}^{m}$ satisfying \eqref{CondSeqOrdSet}, and that
	\eqref{IneqRollGeneral} is  true for the first $m$ intervals $\{I_{k}\}_{k=0}^{m-1}$. From \eqref{IneqRoll},  we obtain
	\begin{align*}
		&\Nzero{Q}{J}+\Ncero{Q}{I_0\setminus J}+\sum_{i=1}^m \Ncero{Q^{(i)}}{I_i}\\
			&\leq\Nzero{Q^{(m-1)}}{J}+\Ncero{Q^{(m-1)}}{\ch{\cup_{i=0}^{m-1}I_i}\setminus J}+m-1+\Ncero{Q^{(m)}}{I_m}\\
&\leq 	\Nzero{Q^{(m)}}{J}+\Ncero{Q^{(m)}}{\intch{\cup_{i=0}^{m-1}I_i}\setminus J}+m+\Ncero{Q^{(m)}}{I_m}\\
			&\leq \Nzero{Q^{(m)}}{J}+\Ncero{Q^{(m)}}{\ch{\cup_{i=0}^{m}I_i}\setminus J}+m.
	\end{align*}
\qed \end{proof}
\begin{corollary}
Under the hypotheses of the above lemma we have
\begin{align}\label{IneqRollzerocero}
\Nzero{Q}{J}+\Ncero{Q}{I_0\setminus J}+\sum_{i=1}^m \Ncero{Q^{(i)}}{I_i}\leq \deg{Q}
\end{align}
for every $J$ closed subinterval of $\intersub{I_0}$. In particular for $J=\emptyset$ we get
\begin{align}\label{IneqRollcero}
\sum_{i=0}^m \Ncero{Q^{(i)}}{I_i}\leq \deg{Q}.
\end{align}
\end{corollary}

\begin{definition}	
	We say that a sequence of ordered pairs $\{(r_i,\nu_i)\}_{i=1}^M \!\subset \RR\times \ZZ_+$ is sequentially-ordered, if  $\nu_1\leq \nu_2 \leq \cdots \leq \nu_M$ and
the set of intervals $ I_{k}=\ch{\{r_i:\nu_i=k\}}$, $k=0,1,\dots, \nu_M, $ satisfy conditions \eqref{CondSeqOrdSet}.
\end{definition}	
	
\begin{lemma}\label{PolMiDeg}
Let $\{(r_i,\nu_i)\}_{i=1}^M\subset \RR\times\ZZp$ be a sequence of $M$ ordered pairs, then there exists a unique monic polynomial $U_M$ of minimal degree ($\leq M$), such that
	\begin{align}\label{cond}
		U_M^{(\nu_i)}(r_i)=0, \quad i=1,2,\dots,M.
	\end{align}
	Furthermore, if $\{(r_i,\nu_i)\}_{i=1}^M$ is sequentially-ordered, then the degree of $U_M$ is  $\MinDeg=\min \mathfrak{I}_{M}-1$, where $$\mathfrak{I}_{M}=\{i: 1\leq i\leq M \text{ and } \nu_i\geq i\}\cup \{M+1\}.$$
\end{lemma}

\begin{proof}
	The existence of a  non-identically-zero polynomial with degree $\leq M$ satisfying \eqref{cond} reduces to solving a homogeneous linear 	system of $M$ equations with $M+1$ unknowns (its coefficients). Thus, a non trivial solution always exists. In addition, if we suppose that 			there exist two different minimal monic polynomials $U_M$ and $\widetilde{U}_M$, then  the polynomial $\widehat{U}_M=U_M-\widetilde{U}_M$ is 	not identically zero, it satisfies \eqref{cond}, and  $\deg{\widehat{U}_M}<\deg{U_M}$. So, if we divide $\widehat{U}_M$ by its leading 				coefficient, we reach a contradiction.

	The rest of the proof runs  by induction on the number of points $M$.  For $M=1$, the result follows taking
		$$U_1(x)=\funD{x-r_1}{\nu_1=0}{1}{\nu_1\geq 1.}$$

	Suppose that,  for each sequentially-ordered sequence  of $M-1$ ordered pairs,  the corresponding minimal  polynomial $U_{M-1}$ has degree $\MinDeg[M-1]$.

	Let $\{(r_i,\nu_i)\}_{i=1}^{M}$ be a sequentially-ordered sequence of $M$ ordered pairs. Obviously,  $\{(r_i,\nu_i)\}_{i=1}^{M-1}$ is a 		sequence of $M-1$ ordered pairs which is sequentially-ordered, $\deg{U_{M}}\geq \deg{U_{M-1}}$, and from the induction hypothesis $\deg{U_{M-1}}=\MinDeg[M-1]$.  Now, we shall split the proof in two cases:
\begin{enumerate}
	\item If $\MinDeg[M]=M$, then for all $1\leq i\leq M$ we have $ \nu_i< i$, which yields
		\begin{equation*}
			\deg{U_{M}}\geq \deg{U_{M-1}}=\MinDeg[M-1]=M-1 \geq \nu_{M}.
		\end{equation*}
		Since $\{(r_i,\nu_i)\}_{i=1}^{M}$ is sequentially-ordered, from \eqref{IneqRollcero} we get
			$$ M\leq \sum_{i=0}^{\nu_{M}} \Ncero{U_{M}^{(i)}}{I_i} \leq  \deg{U_{M}},$$
		which implies that  $\deg{U_{M}}=M=\MinDeg$.
	\item  If $\MinDeg\leq M-1$, then  there exists a minimal $j$ ($1\leq j\leq M$), such that $\nu_j\geq j$,  and $ \nu_i< i$ for all 				$1\leq i\leq j-1$. Therefore,  $\MinDeg=j-1=\MinDeg[M-1]$. From the induction hypothesis
			$$\deg{U_{M-1}}=\MinDeg[M-1]=j-1\leq \nu_j-1\leq \nu_{M}-1,$$
		which gives  $U^{(\nu_{M})}_{M-1}\equiv 0$. Hence,  $U_{M}\equiv U_{M-1}$ and, consequently, we get
			$$\deg{U_{M}}=\deg{U_{M-1}}=\MinDeg[M-1]=\MinDeg.$$
\end{enumerate}
\qed \end{proof}

Note that, in Lemma \ref{PolMiDeg}, the assumption of $\{(r_i,\nu_i)\}_{i=1}^M$ being sequentially ordered is necessary for asserting that the polynomial $U_M$ has degree $\MinDeg$. In fact, if we consider $\{(-1,0),(1,0),(0,1)\}$, which is no sequentially ordered, we get $U_3=x^2-1$ and $\MinDeg[3]=3\neq \deg{U_3}$.

\begin{proof}[Proof of Theorem \ref{ThZeroLocDiscrete}]
	Let  $\xi_1<\xi_2<\cdots <\xi_{\eta}$ be the points on $\inter{\Delta}$ where $P_n$ changes sign and suppose that $\eta<n-\mathbf{d}^*$. Since $\|\cdot\|_{p,\vec{\mu}}$ is sequentially-ordered,  the sequence of $\mathbf{d}^*+\eta$ ordered pairs
		$$\{(r_i,\nu_i)\}_{i=1}^{\mathbf{d}^*+\eta}=\{(\xi_i,0)\}_{i=1}^{\eta}\cup\{(c_j,k): A_{j,k}>0, \ j=1,\dots,N, k=0,\dots,m_j\}$$	
	is sequentially ordered (we can assume without loss of generality that $\nu_1\leq\nu_2\leq \cdots\leq \nu_{\mathbf{d}^*+\eta}$). Consequently, from Lemma \ref{PolMiDeg}, there exists a unique monic polynomial $U_{\mathbf{d}^*+\eta}$ of minimal degree,  such that
	\begin{align}\label{CondNull}
		U_{\mathbf{d}^*+\eta}(\xi_i)&=0,\qquad\text{for } i=1,\dots, \eta;\nonumber\\
		U_{\mathbf{d}^*+\eta}^{(k)}(c_j)&=0,\qquad\text{for } (j,k): A_{j,k}>0;
	\end{align}
	and $ \dsty 		\deg{U_{\mathbf{d}^*+\eta}}=\min \mathfrak{I}_{\mathbf{d}^*+\eta}-1\leq \mathbf{d}^*+\eta,$ 	where
	\begin{equation}\label{DegQ}
		\mathfrak{I}_{\mathbf{d}^*+\eta}=\{i: 1\leq i\leq \mathbf{d}^*+\eta \text{ and } \nu_i \geq i\}\cup \{\mathbf{d}^*+\eta+1\}.
	\end{equation}
	Now, we need to consider the following 	two cases.
	\begin{enumerate}
		\item If $\deg{U_{\mathbf{d}^*+\eta}}=\mathbf{d}^*+\eta$, from \eqref{DegQ},  we get $\deg{U_{\mathbf{d}^*+\eta}}=\mathbf{d}^*+\eta\geq\nu_{\eta+\mathbf{d}^*}+1$. 					Thus, taking $I_i=\Delta_i$, $i=0,1,\dots,m$ and the closed interval $J=[\xi_1,\xi_\eta]\subset \inter{\Delta}\subset \intersub{\Delta_0} $ in \eqref{IneqRollzerocero},  we get
			\begin{align*}
				\mathbf{d}^*+\eta\leq&\sum_{k=0}^{\nu_{\mathbf{d}^*+\eta}} \Ncero{U_{\mathbf{d}^*+\eta}^{(k)}}{\Delta_k} \leq \Nzero{U_{\mathbf{d}^*+\eta}}{J}+\Ncero{U_{\mathbf{d}^*+\eta}}{\Delta_0\setminus J}\\
				  &+\sum_{k=1}^{\nu_{\mathbf{d}^*+\eta}} \Ncero{U_{\mathbf{d}^*+\eta}^{(k)}}{\Delta_k}\leq   \deg{U_{\mathbf{d}^*+\eta}}=\mathbf{d}^*+\eta.
			\end{align*}		
		\item If $\deg{U_{\mathbf{d}^*+\eta}}<\mathbf{d}^*+\eta$,  from \eqref{DegQ},  there exists $1\leq j\leq \mathbf{d}^*+\eta$ such that $\deg{U_{\mathbf{d}^*+\eta}}=j-1$, $\nu_{j}\geq j$ and $\nu_i\leq i-1$ for $i=1,2,\dots,j-1$. Hence, $$\nu_{j-1}+1\leq j-1=\deg{U_{\mathbf{d}^*+\eta}}$$  and, again, from \eqref{IneqRollzerocero} we have
			\begin{align*}
				j-1\leq&\sum_{k=0}^{\nu_{j-1}} \Ncero{U_{\mathbf{d}^*+\eta}^{(k)}}{\Delta_k}\leq \Nzero{U_{\mathbf{d}^*+\eta}}{J}+\Ncero{U_{\mathbf{d}^*+\eta}}{\Delta_0\setminus J}\\
				  &+\sum_{k=1}^{\nu_{j-1}} \Ncero{U_{\mathbf{d}^*+\eta}^{(k)}}{\Delta_k} \leq \deg{U_{\mathbf{d}^*+\eta}}=j-1.
			\end{align*}
		\end{enumerate}
	In both cases, we obtain that $U_{\mathbf{d}^*+\eta}$ has no other zeros in $\Delta_0$ than those given by construction and from $\Ncero{U_{\mathbf{d}^*+\eta}}{J}=\Nzero{U_{\mathbf{d}^*+\eta}}{J}$ we obtain that all the zeros on $\inter{\Delta}$ are simple. Thus, in addition to \eqref{CondNull}, we get that $P_nU_{\mathbf{d}^*+\eta}$ does not change sign on $\inter{\Delta}$. So we have
	\begin{align*}
\langle P_n, U_{\mathbf{d}^*+\eta}\rangle_{p,\mu}=&\int_{\Delta} U_{\mathbf{d}^*+\eta} \, \sgn{P_n}|P_n|^{p-1}  d\mu \\
& +\sum_{j=1}^{N}\sum_{k=0}^{m_j}A_{j,k} U_{\mathbf{d}^*+\eta}^{(k)}(c_{j}) \, \sgn{P_n^{(k)}(c_{j})}|P_n^{(k)}(c_{j})|^{p-1}\\
 			=&\int_{\Delta} U_{\mathbf{d}^*+\eta} \,\sgn{P_n}|P_n|^{p-1}  d\mu\neq0.
	\end{align*}
	Since $\deg{U_{\mathbf{d}^*+\eta}}\leq \mathbf{d}^*+\eta<n$ we arrive at a contradiction  with Theorem \ref{TheoremPProduct}.
\qed \end{proof}


\end{document}